\documentclass[12pt]{article}
\usepackage{mathrsfs}

\usepackage{amssymb}
\usepackage{amsfonts}
\usepackage{amsmath}
\usepackage{amsthm}
\usepackage[margin=0.80in]{geometry}
\usepackage{enumerate}
\usepackage{amstext}
\usepackage{layout}

\numberwithin{equation}{section}

\newtheorem{theorem}{Theorem}[section]
\newtheorem{corollary}{Corollary}[section]
\newtheorem{lemma}{Lemma}[section]

\newtheorem{remark}{Remark}[section]
\newtheorem{example}{Example}[section]

\newtheorem{definition}{Definition}[section]

\begin{document}

\noindent {\bf\large{A Characterization of Chover-Type Law of Iterated Logarithm}}

\vskip 0.3cm

\noindent {\bf Deli Li\footnote{Deli Li, Department of Mathematical
Sciences, Lakehead University, Thunder Bay, Ontario P7B
5E1, Canada\\ e-mail: dli@lakeheadu.ca} $\cdot$ Pingyan Chen\footnote{Pingyan Chen,
Department of Mathematics, Jinan University,
Guangzhou 510630, China\\ e-mail: tchenpy@jnu.edu.cn}
\footnote{Corresponding author: Pingyan Chen (Telephone:
+86~15918852098)}}

\vskip 0.3cm

\noindent {\bf Abstract} Let $0 < \alpha \leq 2$ and $- \infty < \beta < \infty$.
Let $\{X_{n}; n \geq 1 \}$ be a sequence of independent copies of a real-valued random variable $X$
and set $S_{n} = X_{1} + \cdots + X_{n}, ~n \geq 1$. We say $X$ satisfies the $(\alpha, \beta)$-Chover-type
law of the iterated logarithm (and write $X \in CTLIL(\alpha, \beta)$) if $\limsup_{n \rightarrow \infty}
\left| \frac{S_{n}}{n^{1/\alpha}} \right|^{(\log \log n)^{-1}} = e^{\beta}$ almost surely. This paper is
devoted to a characterization of $X \in CTLIL(\alpha, \beta)$. We obtain sets of necessary and sufficient
conditions for $X \in CTLIL(\alpha, \beta)$ for the five cases: $\alpha = 2$ and $0 < \beta < \infty$, $\alpha = 2$
and $\beta = 0$, $1 < \alpha < 2$ and $-\infty < \beta < \infty$, $\alpha = 1$ and $- \infty < \beta < \infty$,
and $0 < \alpha < 1$ and $-\infty < \beta < \infty$. As for the case where $\alpha = 2$ and $-\infty < \beta < 0$,
it is shown that $X \notin CTLIL(2, \beta)$ for any real-valued random variable $X$. As a special case of our results,
a simple and precise characterization of the classical Chover law of the iterated logarithm
(i.e., $X \in CTLIL(\alpha, 1/\alpha)$) is given; that is, $X \in CTLIL(\alpha, 1/\alpha)$ if and only if
$\inf \left \{b:~ \mathbb{E} \left(\frac{|X|^{\alpha}}{(\log (e \vee |X|))^{b\alpha}} \right) < \infty \right\} = 1/\alpha$
where $\mathbb{E}X = 0$ whenever $1 < \alpha \leq 2$.

~\\


\vskip 0.3cm


\vskip 0.3cm

\noindent {\bf Running Head}: Chover-type law of iterated logarithm

\section{Introduction}

Throughout, $\{X_{n}; ~n \geq 1 \}$ is a sequence
of independent copies of a real-valued random variable $X$.
As usual, the partial sums of independent identically distributed (i.i.d.) random variables
$X_{n}, ~n \geq 1$ will be denoted by $S_{n} = \sum_{i=1}^{n} X_{i}, ~n \geq 1$. Write 
$Lx = \log (e \vee x), ~x \geq 0$.

When $X$ has a symmetric stable distribution with exponent $\alpha \in (0, 2)$, i.e.,
$\mathbb{E}\left(e^{itX} \right) = e^{-|t|^{\alpha}}$ for $t \in (-\infty, \infty)$,
Chover proved that
\begin{equation}
\limsup_{n \rightarrow \infty}
\left| \frac{S_{n}}{n^{1/\alpha}} \right|^{(\log \log n)^{-1}} = e^{1/\alpha} ~~\mbox{almost surely (a.s.)}
\end{equation}
This is what we call the classical Chover law of iterated
logarithm (LIL). Since then, several papers have been devoted to
develop the classical Chover LIL. See, for example, Hedye
\cite{Heyde} showed that (1.1) holds when $X$ is in the domain of
normal attraction of a nonnormal stable law, Pakshirajan and
Vasudeva \cite{PV} discussed the limit points of the sequence
$\{|S_{n}/n^{1/\alpha}|^{(\log \log n)^{-1}}; n\geq 2\}$, Kuelbs
and Kurtz \cite{KK} obtained the classical Chover LIL in a Hilbert
space setting, Chen \cite{Chen} obtained the classical Chover LIL
for the weighed sums, Vasudeva \cite{Vasudeva}, Qi and Cheng
\cite{QC}, Peng and Qi \cite{Peng} established the Chover LIL when
$X$ is in the domain of attraction of a nonnormal stable law,
Scheffler \cite{Sch} studied the classical Chover LIL when $X$ is
in the generalized domain of operator semistable attraction of
some nonnormal law, Chen and Hu \cite{CH} extended the results of
Kuelbs and Kurtz \cite{KK} to an arbitrary real separable Banach
space, and so on. It should be to point out that the previous
papers only gave sufficient conditions for the classical Chover
LIL.

Motivated by the previous study of the classical Chover LIL, we introduce a general Chover-type LIL as follows.

\vskip 0.3cm

\begin{definition}
Let $0 < \alpha \leq 2$ and $-\infty < \beta < \infty$. Let $\{X, X_{n}; ~n \geq 1 \}$ be a sequence of
real-valued i.i.d. random variables. We say $X$ satisfies the $(\alpha, \beta)$-Chover-type
law of the iterated logarithm (and write $X \in CTLIL(\alpha, \beta)$) if
\begin{equation}
\limsup_{n \rightarrow \infty}
\left| \frac{S_{n}}{n^{1/\alpha}} \right|^{(\log \log n)^{-1}} = e^{\beta} ~~\mbox{a.s.}
\end{equation}
\end{definition}

\vskip 0.2cm

From the classical Chover LIL and Definition 1.1, we see that $X \in CTLIL(\alpha, 1/\alpha)$
(i.e., (1.2) holds with $\beta = 1/\alpha$) when $X$ has a symmetric stable distribution with
exponent $\alpha \in (0, 2)$.

This paper is devoted to a characterization of $X \in
CTLIL(\alpha, \beta)$. The main results are stated in Section 2.
We obtain sets of necessary and sufficient conditions for $X \in
CTLIL(\alpha, \beta)$ for the five cases: $\alpha = 2$ and $0 <
\beta < \infty$ ( see Theorem 2.1), $\alpha = 2$ and $\beta = 0$
(see Theorem 2.2), $1 < \alpha < 2$ and $-\infty < \beta < \infty$
(see Theorem 2.3), $\alpha = 1$ and $- \infty < \beta < \infty$
(see Theorem 2.4), and $0 < \alpha < 1$ and $-\infty < \beta <
\infty$ (see Theorem 2.5). The proofs of Theorems 2.1-2.5 are
given in Section 4. For proving Theorems 2.1-2.5, three
preliminary lemmas are stated in Section 3. Some llustrative
examples are provided in Section 5.

\section{Statement of the Main Results}

The main results of this paper are the following five theorems. We begin with the case where $\alpha = 2$
and $0 < \beta < \infty$.

\begin{theorem}
Let $ 0 < \beta < \infty$. Let $\{X, X_{n}; n \geq 1 \}$ be a sequence of i.i.d. real-valued random variables.
Then
\begin{equation}
X \in CTLIL(2, \beta), ~\mbox{i.e.,}~ \limsup_{n \rightarrow \infty} \left| \frac{S_{n}}{\sqrt{n}} \right|^{(\log \log n)^{-1}}
= e^{\beta} ~~\mbox{a.s.}
\end{equation}
if and only if
\begin{equation}
\mathbb{E}X = 0~~\mbox{and}~~\inf \left \{b > 0:~ \mathbb{E} \left(\frac{X^{2}}{(L|X|)^{2b}} \right) < \infty \right\} = \beta.
\end{equation}
\end{theorem}

\vskip 0.3cm

For the case where $\alpha = 2$ and $\beta = 0$, we have the following result.

\begin{theorem}
Let $\{X, X_{n}; n \geq 1 \}$ be a sequence of i.i.d. non-degenerate real-valued random variables.
Then
\begin{equation}
\limsup_{n \rightarrow \infty} \left| \frac{S_{n}}{\sqrt{n}} \right|^{(\log \log n)^{-1}} \leq 1 ~~\mbox{a.s.}
\end{equation}
if and only if
\begin{equation}
\mathbb{E}X = 0~~\mbox{and}~~\inf \left \{b > 0:~ \mathbb{E} \left(\frac{X^{2}}{(L|X|)^{2b}} \right) < \infty \right\} = 0.
\end{equation}
In either case, we have
\begin{equation}
X \in CTLIL(2, 0), ~\mbox{i.e.,}~
\limsup_{n \rightarrow \infty} \left| \frac{S_{n}}{\sqrt{n}} \right|^{(\log \log n)^{-1}} = 1 ~~\mbox{a.s.}
\end{equation}
\end{theorem}

\vskip 0.3cm

\begin{remark}
Let $c$ be a constant. Note that
\[
\limsup_{n \rightarrow \infty} \left| \frac{nc}{\sqrt{n}} \right|^{(\log \log n)^{-1}}
= \left \{
\begin{array}{ll}
0 & if\ c = 0,\\
& \\
\infty & if\ c \neq 0.
\end{array}
\right.
\]
Thus, from Theorem 2.2, we conclude that, for any $- \infty < \beta < 0$, $X \not\in CTLIL(2, \beta)$ for
any real-valued random variable $X$.
\end{remark}

\vskip 0.3cm

In the next three theorems, we provide necessary and sufficient conditions for $X \in CTLIL(\alpha, \beta)$
for the three cases where $1 < \alpha < 2$ and $-\infty < \beta < \infty$, $\alpha = 1$ and $- \infty < \beta < \infty$,
and $0 < \alpha < 1$ and $-\infty < \beta < \infty$ respectively.

\begin{theorem}
Let $1 < \alpha < 2$ and $ - \infty < \beta < \infty$.
Let $\{X, X_{n}; n \geq 1 \}$ be a sequence of i.i.d. real-valued random variables. Then
\[
X \in CTLIL(\alpha, \beta), ~\mbox{i.e.,}~
\limsup_{n \rightarrow \infty} \left| \frac{S_{n}}{n^{1/\alpha}} \right|^{(\log \log n)^{-1}} = e^{\beta} ~~\mbox{a.s.}
\]
if and only if
\[
\mathbb{E}X = 0~~\mbox{and}~~\inf \left \{b:~ \mathbb{E} \left(\frac{|X|^{\alpha}}{(L|X|)^{b\alpha}} \right) < \infty \right\} = \beta.
\]
\end{theorem}

\vskip 0.3cm

\begin{theorem}
Let $ - \infty < \beta < \infty$.
Let $\{X, X_{n}; n \geq 1 \}$ be a sequence of i.i.d. real-valued random variables. Then
\[
X \in CTLIL(1, \beta), ~\mbox{i.e.,}~
\limsup_{n \rightarrow \infty} \left| \frac{S_{n}}{n} \right|^{(\log \log n)^{-1}} = e^{\beta} ~~\mbox{a.s.}
\]
if and only if
\[
\left \{
\begin{array}{ll}
\mbox{$\displaystyle
\inf \left \{b:~ \mathbb{E} \left(\frac{|X|}{(L|X|)^{b}} \right) < \infty \right\} = \beta$} & \mbox{if $\beta > 0$,}\\
&\\
\mbox{$\displaystyle \mbox{either}~~\mathbb{E}|X| < \infty ~\mbox{and}~\mathbb{E}X \neq 0
~~\mbox{or}~~\inf \left \{b:~ \mathbb{E} \left(\frac{|X|}{(L|X|)^{b}} \right) < \infty \right\} = 0$} & \mbox{if $\beta = 0$,}\\
&\\
\mbox{$\displaystyle \mathbb{E}X = 0
~\mbox{and}~~\inf \left \{b:~ \mathbb{E} \left(\frac{|X|}{(L|X|)^{b}} \right) < \infty \right\} = \beta$} & \mbox{if $\beta < 0$.}
\end{array}
\right.
\]
In particular, $\mathbb{E}|X|<\infty$ and $\mathbb{E}X \neq 0$
imply that
\[
\lim_{n\rightarrow\infty}\left|\frac{S_n}{n}\right|^{(\log\log n)^{-1}}=1\ \ a.s.
\]
\end{theorem}

\vskip 0.3cm

\begin{theorem}
Let $0 < \alpha < 1$ and $ - \infty < \beta < \infty$.
Let $\{X, X_{n}; n \geq 1 \}$ be a sequence of i.i.d. real-valued random variables. Then
\[
X \in CTLIL(\alpha, \beta),~\mbox{i.e.,}~
\limsup_{n \rightarrow \infty} \left| \frac{S_{n}}{n^{1/\alpha}} \right|^{(\log \log n)^{-1}} = e^{\beta} ~~\mbox{a.s.}
\]
if and only if
\[
\inf \left \{b:~ \mathbb{E} \left(\frac{|X|^{\alpha}}{(L|X|)^{b\alpha}} \right) < \infty \right\} = \beta.
\]
\end{theorem}

\vskip 0.3cm

\begin{remark}
From our Theorems 2.1, 2.3, 2.4, and 2.5, a simple and precise characterization of the classical Chover LIL
(i.e., $X \in CTLIL(\alpha, 1/\alpha)$) is obtained as follows. For $0 < \alpha \leq 2$, we have
\[
X \in CTLIL(\alpha, 1/\alpha), ~\mbox{i.e.,}~
\limsup_{n \rightarrow \infty} \left| \frac{S_{n}}{n^{1/\alpha}} \right|^{(\log \log n)^{-1}} = e^{1/\alpha} ~~\mbox{a.s.}
\]
if and only if
\[
\inf \left \{b:~ \mathbb{E} \left(\frac{|X|^{\alpha}}{(L|X|)^{b\alpha}} \right) < \infty \right\} = 1/\alpha
~~\mbox{where}~~\mathbb{E}X = 0~~\mbox{whenever}~~1 < \alpha \leq 2.
\]
\end{remark}

\vskip 0.3cm

Our Theorems 2.1-2.5 also imply the following two interesting results.

\vskip 0.2cm

\begin{corollary}
Let $0 < \alpha \leq 2$. Let $\{X, X_{n}; n \geq 1 \}$ be a sequence of i.i.d. real-valued random variables. Then
\[
\lim_{n \rightarrow \infty} \left| \frac{S_{n}}{n^{1/\alpha}} \right|^{(\log \log n)^{-1}} = 0 ~~\mbox{a.s.}
\]
if and only if
\[
\left \{
\begin{array}{ll}
\mbox{$\displaystyle \inf \left \{b:~ \mathbb{E} \left(\frac{|X|^{\alpha}}{(L|X|)^{b\alpha}} \right) 
< \infty \right\} =  - \infty $} & \mbox{if $\displaystyle 0 < \alpha < 1$}, \\
& \\
\mbox{$\displaystyle \mathbb{E}X = 0 ~\mbox{and}~ 
\inf \left \{b:~ \mathbb{E} \left(\frac{|X|^{\alpha}}{(L|X|)^{b\alpha}} \right) < \infty \right\} = - \infty $}
& \mbox{if $\displaystyle 1 \leq \alpha < 2$}, \\
& \\
\mbox{$\displaystyle X = 0$ ~a.s.} & \mbox{if $\displaystyle \alpha = 2$}.
\end{array}
\right.
\]
\end{corollary}

\vskip 0.2cm

\begin{corollary}
Let $0 < \alpha \leq 2$. Let $\{X, X_{n}; n \geq 1 \}$ be a sequence of i.i.d. real-valued random variables. Then
\[
\limsup_{n \rightarrow \infty} \left| \frac{S_{n}}{n^{1/\alpha}} \right|^{(\log \log n)^{-1}} = \infty~~\mbox{a.s.}
\]
if and only if
\[
\left \{
\begin{array}{ll}
\mbox{$\displaystyle \inf \left \{b:~ \mathbb{E} \left(\frac{|X|^{\alpha}}{(L|X|)^{b\alpha}} \right) 
< \infty \right\} =  \infty $} & \mbox{if $\displaystyle 0 < \alpha \leq 1$},\\
& \\
\mbox{eithor $\displaystyle \mathbb{E}X \neq 0$ or}~ \mbox{$\displaystyle 
\inf \left \{b:~ \mathbb{E} \left(\frac{|X|^{\alpha}}{(L|X|)^{b\alpha}} \right) < \infty \right\} = \infty$}
& \mbox{if $\displaystyle 1 < \alpha \leq 2$}.
\end{array}
\right.
\]
\end{corollary}

\section{Preliminary lemmas}

To prove the main results, we use the following three preliminary lemmas.
The first lemma is new and may be of independent interest.

\vskip 0.3cm

\begin{lemma}
Let $\{a_{n};~n \geq 1 \}$ be a sequence of real numbers. Let $\{c_{n};~n \geq 1 \}$ be a sequence
of positive real numbers such that
\begin{equation}
\lim_{n \rightarrow \infty} c_{n} = \infty.
\end{equation}
Then we have

(i)~~There exists a constant $-\infty < \beta < \infty$ such that
\begin{equation}
\limsup_{n \rightarrow \infty} \left|a_{n} \right|^{1/c_{n}} = e^{\beta}
\end{equation}
if and only if
\begin{equation}
\limsup_{n \rightarrow \infty}
\frac{\left|a_{n} \right|}{e^{bc_{n}}} =
\left \{
\begin{array}{ll}
0 & \mbox{for all}~ b > \beta,\\
&\\
\infty & \mbox{for all}~ b < \beta;
\end{array}
\right.
\end{equation}

(ii)~~There exists a constant $-\infty < \beta < \infty$ such that
\[
\liminf_{n \rightarrow \infty} \left|a_{n} \right|^{1/c_{n}} = e^{\beta}
\]
if and only if
\[
\liminf_{n \rightarrow \infty}
\frac{\left|a_{n} \right|}{e^{bc_{n}}} =
\left \{
\begin{array}{ll}
0 & \mbox{for all}~ b > \beta,\\
&\\
\infty & \mbox{for all}~ b < \beta;
\end{array}
\right.
\]

(iii)~~There exists a constant $-\infty < \beta < \infty$ such that
\[
\lim_{n \rightarrow \infty} \left|a_{n} \right|^{1/c_{n}} = e^{\beta}
\]
if and only if
\[
\lim_{n \rightarrow \infty}
\frac{\left|a_{n} \right|}{e^{bc_{n}}} =
\left \{
\begin{array}{ll}
0 & \mbox{for all}~ b > \beta,\\
&\\
\infty & \mbox{for all}~ b < \beta.
\end{array}
\right.
\]
\end{lemma}

{\it Proof}~~We prove the sufficiency part of Part (i) first. It follows from (3.3) that the set
\[
\left\{n \geq 1;~\frac{\left|a_{n} \right|}{e^{bc_{n}}} > 1 \right \}
\left \{
\begin{array}{ll}
\mbox{has finitely many elements}~ &\mbox{if}~b > \beta,\\
&\\
\mbox{has infinitely many elements}~ &\mbox{if}~b < \beta.
\end{array}
\right.
\]
Note that
\[
\left\{n \geq 1;~\left|a_{n} \right|^{1/c_{n}} > e^{b} \right \}
= \left\{n \geq 1;~\frac{\left|a_{n} \right|}{e^{bc_{n}}} > 1 \right \}.
\]
We thus conclude that
\[
\limsup_{n \rightarrow \infty} \left|a_{n} \right|^{1/c_{n}}
\left\{
\begin{array}{ll}
\mbox{$\displaystyle \leq e^{b}$}~ & \mbox{for all}~ b > \beta,\\
&\\
\mbox{$\displaystyle \geq e^{b}$}~ & \mbox{for all}~ b < \beta
\end{array}
\right.
\]
which ensures (3.2).

We now prove the necessity part of Part (i). For all $b \neq \beta$, let $h = (b + \beta)/2$. Then
\[
\left\{
\begin{array}{ll}
\beta < h < b & \mbox{if}~b > \beta,\\
&\\
b < h < \beta & \mbox{if}~b < \beta.
\end{array}
\right.
\]
It follows from (3.2) that the set
\[
\left\{n \geq 1;~\left|a_{n} \right|^{1/c_{n}} > e^{h} \right \}
\left \{
\begin{array}{ll}
\mbox{has finitely many elements}~ &\mbox{if}~b > \beta,\\
&\\
\mbox{has infinitely many elements}~ &\mbox{if}~b < \beta.
\end{array}
\right.
\]
Note that
\[
\left\{n \geq 1;~\frac{\left|a_{n} \right|}{e^{hc_{n}}} > 1 \right \}
= \left\{n \geq 1;~\left|a_{n} \right|^{1/c_{n}} > e^{h} \right \}.
\]
We thus have that
\begin{equation}
\limsup_{n \rightarrow \infty} \frac{\left|a_{n} \right|}{e^{hc_{n}}}
\left\{
\begin{array}{ll}
\leq 1~ & \mbox{if}~ b > \beta,\\
&\\
\geq 1~ & \mbox{if}~ b < \beta.
\end{array}
\right.
\end{equation}
Note that (3.1) implies that
\[
\lim_{n \rightarrow \infty} \frac{e^{hc_{n}}}{e^{bc_{n}}} = \lim_{n \rightarrow \infty} e^{(h-b)c_{n}}
= \left \{
\begin{array}{ll}
0 & \mbox{if}~b > \beta,\\
&\\
\infty & \mbox{if}~ b < \beta.
\end{array}
\right.
\]
Thus it follows from (3.4) that
\[
\limsup_{n \rightarrow \infty} \frac{\left|a_{n} \right|}{e^{bc_{n}}}
= \left(\lim_{n \rightarrow \infty} \frac{e^{hc_{n}}}{e^{bc_{n}}} \right)
\left(\limsup_{n \rightarrow \infty} \frac{\left|a_{n} \right|}{e^{hc_{n}}} \right)
= \left\{
\begin{array}{ll}
0~ & \mbox{if}~ b > \beta,\\
&\\
\infty~ & \mbox{if}~ b < \beta,
\end{array}
\right.
\]
i.e., (3.3) holds.

We leave the proofs of Parts (ii) and (iii) to the reader since they are similar to the proof of Part (i).
This completes the proof of Lemma 3.1. ~$\Box$

\vskip 0.3cm

The following result is a special case of Corollary 2 of Einmahl
and Li \cite{EL}.

\vskip 0.3cm

\begin{lemma}
Let $b > 0$. Let $\{X, X_{n}; n \geq 1 \}$ be a sequence of i.i.d. random variables.
Then
\[
\lim_{n \rightarrow \infty} \frac{S_{n}}{\sqrt{n} (\log n)^{b}} = 0 ~~\mbox{a.s.}
\]
if and only if
\[
\mathbb{E}X = 0, ~~\mathbb{E}\left(\frac{X^{2}}{(L|X|)^{2b}} \right) < \infty, ~~
\mbox{and}~ \lim_{x \rightarrow \infty} \frac{LL x}{(L x)^{2b}} H(x) = 0,
\]
where $H(x) = \mathbb{E}\left(X^{2}I\{|X| \leq x \} \right)$, $x \geq 0$.
\end{lemma}

\vskip 0.3cm

The following result is a generalization of
Kolmogorov-Marcinkiewicz-Zygmund strong law of large numbers and
follows easily from Theorems 1 and 2 of Feller \cite{Feller}.

\vskip 0.3cm

\begin{lemma}
Let $0 < \alpha < 2$ and $- \infty < b < \infty$. Let $\{X, X_{n}; n \geq 1 \}$ be a sequence of i.i.d. random variables.
Then
\[
\lim_{n \rightarrow \infty} \frac{S_{n}}{n^{1/\alpha} (\log n)^{b}} = 0 ~~\mbox{a.s.}
\]
if and only if
\[
\mathbb{E}\left(\frac{|X|^{\alpha}}{(L|X|)^{b\alpha}} \right) < \infty
\]
where $\mathbb{E}X = 0$ whenever either $1 < \alpha < 2$ or $\alpha = 1$ and $- \infty < b \leq 0$.
\end{lemma}

\section{Proofs of the Main Results}

In this section, we only give the proofs of Theorems 2.1-2.2. By applying Lemmas 3.1 and 3.3,
the proofs of Theorems 2.3-2.5 involves only minor modifications of the proof of Theorem 2.1
and will be omitted.

\vskip 0.2cm

{\it Proof of Theorem 2.1}~~We prove the sufficiency part first. Note that the second part of (2.2) implies that
\begin{equation}
\mathbb{E}\left(\frac{X^{2}}{(L|X|)^{2b}} \right) < \infty~~\mbox{for all}~b > \beta.
\end{equation}
We thus see that for all $b > \beta$,
\[
\begin{array}{lll}
\mbox{$\displaystyle H(x)$}
& = & \mbox{$\displaystyle \mathbb{E}\left(X^{2}I\{|X| \leq x \} \right)$}\\
&&\\
& \leq & \mbox{$\displaystyle \mathbb{E} \left(\frac{X^{2}}{(L|X|)^{2h}} (Lx)^{2h} I\{|X| \leq x \} \right)$}\\
&&\\
& \leq & \mbox{$\displaystyle (Lx)^{2h} \mathbb{E} \left(\frac{X^{2}}{(L|X|)^{2h}} \right)$},
\end{array}
\]
where $h = (b + \beta)/2$. Since $\beta < h < b$, it follows from (4.1) that
\[
\frac{LL x}{(L x)^{2b}} H(x) \leq \frac{LLx}{(Lx)^{2(b-h)}} \mathbb{E} \left(\frac{X^{2}}{(L|X|)^{2h}} \right)
\rightarrow 0 ~~\mbox{as}~ x \rightarrow \infty.
\]
We thus conclude from (2.2) that
\[
\mathbb{E}X = 0, ~~\mathbb{E}\left(\frac{X^{2}}{(L|X|)^{2b}} \right) < \infty,
~~\mbox{and}~ \lim_{x \rightarrow \infty} \frac{LL x}{(L x)^{2b}} H(x) = 0 ~~\mbox{for all}~b > \beta > 0
\]
which, by applying Lemma 3.2, ensures that
\begin{equation}
\lim_{n \rightarrow \infty} \frac{S_{n}}{\sqrt{n} (\log n)^{b}} = 0 ~~\mbox{a.s. ~ for all}~b > \beta.
\end{equation}
Since $\beta > 0$, the second part of (2.2) implies that
\[
\mathbb{E}\left(\frac{X^{2}}{(L|X|)^{2b}} \right) = \infty~~\mbox{for all}~b < \beta
\]
which ensures that
\begin{equation}
\limsup_{n \rightarrow \infty} \frac{\left|S_{n} \right|}{\sqrt{n} (\log n)^{b}}
= \infty ~~\mbox{a.s. ~ for all}~b < \beta.
\end{equation}
Let $A_{n} = S_{n}/\sqrt{n}$ and $c_{n} = LLn$, $n \geq 1$. It then follows from (4.2) and (4.3) that
\begin{equation}
\limsup_{n \rightarrow \infty} \frac{\left|A_{n} \right|}{e^{bc_{n}}}
= \limsup_{n \rightarrow \infty} \frac{\left|S_{n} \right|}{\sqrt{n} (\log n)^{b}}
= \left \{
\begin{array}{l}
0 ~~\mbox{a.s. for all}~ b > \beta,\\
\\
\infty ~~\mbox{a.s. for all}~ b < \beta.
\end{array}
\right.
\end{equation}
By Lemma 3.1, we see that (4.4) is equivalent to
\[
\limsup_{n \rightarrow \infty} \left|A_{n}\right|^{1/c_{n}} = e^{\beta}~~\mbox{a.s.},
\]
i.e., (2.1) holds.

We now prove the necessity part. By Lemma 3.1, (2.1) is equivalent to (4.4) which ensures that
(4.2) holds. By Lemma 3.2, we conclude from (4.2) that
\begin{equation}
\mathbb{E}X = 0 ~~\mbox{and}~~\mathbb{E}\left(\frac{X^{2}}{(L|X|)^{2b}} \right) < \infty
~~\mbox{for all}~b > \beta.
\end{equation}
Since $0 < \beta < \infty$, it follows from (4.5) that
\[
\beta_{1} \stackrel{\Delta}{=}
\inf\left\{b > 0: ~\mathbb{E}\left(\frac{X^{2}}{(L|X|)^{2b}} \right) < \infty \right \} \leq \beta.
\]
If $\beta_{1} < \beta$ then, using the argument in the proof of the sufficiency part, we have that
\begin{equation}
\lim_{n \rightarrow \infty} \frac{S_{n}}{\sqrt{n} (\log n)^{b}} = 0 ~~\mbox{a.s. ~ for all}~b > \beta_{1}.
\end{equation}
Hence, by Lemma 3.1, (4.6) implies that
\[
\limsup_{n \rightarrow \infty} \left|\frac{S_{n}}{\sqrt{n}} \right|^{(\log \log n)^{-1}} \leq e^{\beta_{1}}
< e^{\beta} ~~\mbox{a.s.}
\]
which is in contradiction to (2.1). Thus (2.2) holds. The proof of Theorem 2.1 is complete. ~$\Box$

\vskip 0.3cm

{\it Proof of Theorem 2.2}~~Using the same argument used in the proof of the sufficiency part
of Theorem 2.1, we have from (2.4) that
\begin{equation}
\lim_{n \rightarrow \infty} \frac{S_{n}}{\sqrt{n} (\log n)^{b}} = 0 ~~\mbox{a.s. ~ for all}~b > 0.
\end{equation}
Since $X$ is a non-degenerate random variable, by the classical Hartman-Wintner-Strassen
LIL, we have that
\[
\limsup_{n \rightarrow \infty} \frac{|S_{n} |}{\sqrt{2n \log \log n}} > 0 ~~\mbox{a.s.}
\]
which implies that
\begin{equation}
\limsup_{n \rightarrow \infty} \frac{|S_{n}|}{\sqrt{n} (\log n)^{b}} = \infty ~~\mbox{a.s. ~ for all}~b < 0.
\end{equation}
Let $A_{n} = S_{n}/\sqrt{n}$ and $c_{n} = LLn$, $n \geq 1$. It then follows from (4.7) and (4.8) that
\begin{equation}
\limsup_{n \rightarrow \infty} \frac{\left|A_{n} \right|}{e^{bc_{n}}}
= \limsup_{n \rightarrow \infty} \frac{\left|S_{n} \right|}{\sqrt{n} (\log n)^{b}}
= \left \{
\begin{array}{l}
0 ~~\mbox{a.s. for all}~ b > 0,\\
\\
\infty ~~\mbox{a.s. for all}~ b < 0.
\end{array}
\right.
\end{equation}
By Lemma 3.1, we see that (4.9) is equivalent to
\[
\limsup_{n \rightarrow \infty} \left|A_{n}\right|^{1/c_{n}} = e^{0} = 1~~\mbox{a.s.},
\]
i.e., (2.5) holds, so does (2.3).

Using the same argument used in the proof of the necessity part of Theorem 2.1, we conclude from (2.3) that
\[
\mathbb{E}X = 0 ~~\mbox{and}~~\inf\left\{b > 0: ~\mathbb{E}\left(\frac{X^{2}}{(L|X|)^{2b}} \right) < \infty \right \} \leq 0.
\]
Clearly
\[
\inf\left\{b > 0: ~\mathbb{E}\left(\frac{X^{2}}{(L|X|)^{2b}} \right) < \infty \right \} \geq 0.
\]
Thus (2.4) holds. The proof of Theorem 2.2 is therefore complete. ~$\Box$

\section{Examples}

In this section, we provide the following examples to illustrate our main results.
By applying Theorems 2.3-2.5, we rededuce the classical Chover LIL in the first
example.

\vskip 0.3cm

\begin{example}
Let $0 < \alpha \leq 2$. Let $X$ be a symmetric real-valued stable random variable
with exponent $\alpha$. Clearly, $\mathbb{E}X = 0$ whenever $1 < \alpha \leq 2$.

For $0 < \alpha < 2$, we have
\[
\mathbb{P}(|X| > x) \sim \left(\frac{\sin\left(\pi \alpha/2 \right) \Gamma(\alpha)}{\pi} \right)
|x|^{-\alpha} ~~\mbox{as}~~x \rightarrow \infty,
\]
it follows that
\[
\mathbb{E} \left(\frac{|X|^{\alpha}}{(L|X|)^{b\alpha}} \right)
\left \{
\begin{array}{ll}
< \infty & if\ b > 1/\alpha,\\
& \\
= \infty & if\ b \leq 1/\alpha
\end{array}
\right.
\]
and hence that
\[
\inf \left \{b:~ \mathbb{E} \left(\frac{|X|^{\alpha}}{(L|X|)^{b\alpha}} \right) < \infty \right\}
= 1/\alpha.
\]
Thus, by Theorems 2.3-2.5, $X \in CTLIL(\alpha, 1/\alpha)$ (i.e., the classical
Chover LIL follows).

However, for $\alpha = 2$, we have $\mathbb{E}X^{2} = 1$. Hence, by Theorems 2.1 and 2.2, we
see that $X \notin CTLIL(2, 1/2)$ but $X \in CTLIL(2, 0)$.
\end{example}

\vskip 0.3cm

From our second example, we will see that $X \in CTLIL(\alpha, \beta)$ for some certain
$\alpha$ and $\beta$ even if the distribution of $X$ is not in the domain of attraction
of the stable distribution with exponent $\alpha$.

\vskip 0.3cm

\begin{example}
Let $0 < \alpha \leq 2$. Let $d_{n} = \exp\left(2^{n}\right), ~n \geq 1$. Given
$- \infty < \lambda < \infty$. Let $X$ be a symmetric i.i.d. real-valued random variable
such that
\[
\mathbb{P}\left(X = -d_{n} \right) = \mathbb{P}\left(X =
d_{n} \right) = \left(\frac{c}{2} \right)
\frac{\log^{\lambda}d_{n}}{d_{n}^{\alpha}} = \left(\frac{c}{2} \right)
\frac{2^{n\lambda}}{d_{n}^{\alpha}}, ~~n \geq 1
\]
where
\[
c = c(\alpha, \lambda) = \left(\sum_{n = 1}^{\infty}
\frac{\log^{\lambda}d_{n}}{d_{n}^{\alpha}} \right)^{-1} > 0.
\]
Then
\[
\limsup_{x \to \infty} \frac{x^{\alpha}}{\log^{\lambda}x} \mathbb{P} (|X| \geq x ) = c > 0
~~\mbox{and}~~
\liminf_{x \to \infty} \frac{x^{\alpha}}{\log^{\lambda}x} \mathbb{P} (|X| \geq x ) = 0.
\]
Thus the distribution of $X$ is not in the domain of attraction of the stable distribution
with exponent $\alpha$. Also $\mathbb{E}X = 0$ whenever either $1 < \alpha \leq 2$ or $\alpha = 1$
and $\lambda < 0$. It is easy to see that
\[
\mathbb{E} \left(\frac{|X|^{\alpha}}{(L|X|)^{b\alpha}} \right) = \left(\frac{c}{2} \right)
\sum_{n = 1}^{\infty} \left(2^{\lambda - b \alpha} \right)^{n}
\left \{
\begin{array}{ll}
< \infty & if\ b > \lambda/\alpha,\\
& \\
= \infty & if\ b \leq \lambda/\alpha
\end{array}
\right.
\]
and hence that
\[
\inf \left \{b:~ \mathbb{E} \left(\frac{|X|^{\alpha}}{(L|X|)^{b\alpha}} \right) < \infty \right\}
= \lambda/\alpha.
\]
Thus, by Theorems 2.1-2.5, we have

{\rm (1)} If $\alpha = 2$ and $0 < \lambda < \infty$, then $X \in CTLIL(2, \lambda/\alpha)$.

{\rm (2)} If $\alpha=2$ and $ -\infty < \lambda \leq 0$, then $X \in CTLIL(2, 0)$.

{\rm (3)} If $0 < \alpha <2$, then $X \in CTLIL(\alpha, \lambda/\alpha)$.
\end{example}

\vskip 0.3cm

Our third example shows that $X$ may satisfy the other Chover-type LIL studied by Chen and Hu \cite{CH}
when $X \notin CTLIL(\alpha, \beta)$.

\vskip 0.3cm

\begin{example}
Define the density function $f(x)$ of $X$ by
\[
f(x)= \left \{
\begin{array}{ll}
0 & if\ 0\leq |x|<e,\\
&\\
\displaystyle \frac{c}{|x|^{\alpha+1}}\exp\left(p(\log
|x|)^\gamma\right) & if\ |x|\geq e,
\end{array}
\right.
\]
where $p \neq 0$, $0 < \gamma < 1$, $c = c(\alpha, p, \gamma)$ is
a positive constant such that such that
$\int^\infty_{-\infty}f(x)dx=1$. On simplification one can show
that for any $-\infty < b < \infty$
\[
\displaystyle
\mathbb{E}\left(\frac{|X|^\alpha}{(L|X|)^{b\alpha}}\right)=2c\int^\infty_e\frac{\exp\left(p(\log
x)^\gamma\right)}{x(\log x)^{b\alpha}}dx  \left \{
\begin{array}{ll}
<\infty & if\ p<0,\\
&\\
\displaystyle =\infty & if\ p>0.
\end{array}
\right.
\]
From Theorem 2.1-2.5, we have

{\rm (1)} If $\alpha=2$ and $p < 0$, then $X\in CTLIL(2, 0)$.

{\rm (2)} If $\alpha=2$ and $p > 0$, then $X\notin CTLIL(2, \beta)$ for any
$0 \leq \beta<\infty$.

{\rm (3)} If $0 < \alpha < 2$, then $X \notin CTLIL(\alpha, \beta)$ for any
$- \infty <\beta<\infty$.
\\
However, for $0 < \alpha < 2$, by Theorem 3.1 in Chen and Hu {\rm \cite{CH}}, we have
\[
\limsup_{n\rightarrow\infty}\left|\frac{S_n}{B(n)}\right|^{(\log\log n)^{-1}}=e^{1/\alpha}\ \ a.s.,
\]
where $B(x)$ is the inverse function of $x^\alpha/\exp\left(p(\log
x)^\gamma\right)$, $x\geq e$.
\end{example}

\vskip 0.3cm

\noindent {\bf Acknowledgements}

\noindent The research of Deli Li was partially supported
by a grant from the Natural Sciences and Engineering Research Council of
Canada and the research of Pingyan Chen was partially supported
by the National Natural Science Foundation of China (No. 11271161)

\vskip 0.5cm

\bibliographystyle{plain}

\begin{thebibliography}{10}
\bibitem{Chen} Chen, P.: {\it Limiting behavior of weighted sums with stable
        distributions}. Statist. Probab. Lett.  {\bf 60}, 367-375 (2002).
\bibitem{CH} Chen, P., Hu, T.-C.: {\it Limiting behavior for random elements with heavy
        tail}. Taiwanese J.  Math.  {\bf 16}, 217-236 (2012).
\bibitem{Chover} Chover, J.: {\it A law of the iterated logarithm for stable
        summands}. Proc. Amer. Math. Soc. {\bf 17},441-443 (1966).
\bibitem{EL} Einmahl, U., Li, D.: {\it Some results on two-sided LIL
        behavior}. Ann. Probab. {\bf 33}, 1601-1624 (2005).
\bibitem{Feller} Feller, W.: {\it A Limit Theoerm for Random Variables with Infinite
        Moments}. Amer. J. Math. {\bf 68}, 257-262 (1946).
\bibitem{Heyde} Heyde, C. C.: {\it A note concerning behaviour of iterated
        logarithm type}. Proc. Amer. Math. Soc. {\bf 23}, 85-90 (1969).
\bibitem{KK} Kuelbs, J., Kurtz, T.: {\it Berry-Esseen Estimates in Hilbert Space
        and an Application to the Law of the Iterated Logarithm}. Ann. Probab.
        {\bf 2}, 387-407 (1974).
\bibitem{PV} Pakshirajan, R. P., Vasudeva, R.: {\it A law of the iterated
        logarithm for stable summands}. Tran. Amer. Math. Soc. {\bf 232}, 33-42 (1977).
\bibitem{Peng} Peng, L., Qi, Y.: {\it Chover-type laws of the iterated logarithm
        for weighted sums}. Statist. Probab. Lett. {\bf 65}, 401-410 (2003).
\bibitem{QC} Qi, Y., Cheng, P.: {\it On the law of the iterated logarithm for
        the partial sum in the domain of attraction of stable
        distribution}. Chinese Ann. Math. {\bf 17}(A), 195-206 (1996)(in Chinese).
\bibitem{Sch} Scheffler, H.-P.: {\it A law of the iterated logarithm for
        heavy-tailed random vectors}. Probab. Theory. Relat.
        Fields, {\bf 116}, 257-271 (2000).
\bibitem{Vasudeva} Vasudeva, R.: {\it Chover's law of the iterated logarithm and weak
        convergence}. Acta Math. Hung. {\bf 44}, 215-221 (1984).
\end{thebibliography}

\end{document}